\documentclass{amsart}
\usepackage{amsmath, amssymb, amscd}

\newcommand{\et}{{ {\text {\'{e}}t} }}
\newcommand{\cO}{{\mathcal O}}

\newcommand{\fm}{{\mathfrak m}}
\newcommand{\fp}{{\mathfrak p}}

\newcommand{\fG}{{\mathfrak G}}

\newcommand{\Gm}{{\mathbb G}_m}

\newcommand{\CC}{{\mathbb C}}
\newcommand{\U}{{\mathbb U}}
\newcommand{\FF}{{\mathbb F}}
\newcommand{\NN}{{\mathbb N}}
\newcommand{\ZZ}{{\mathbb Z}}

\newcommand{\alg}{{\rm alg}}
\newcommand{\tor}{{\rm tor}}

\newcommand{\ACFA}{{\rm ACFA}}
\newcommand{\Fix}{{\rm Fix}}

\newtheorem{thm}{Theorem}[section]
\newtheorem{lem}[thm]{Lemma}

\theoremstyle{definition}
\newtheorem{Def}[thm]{Definition}

\begin{document}

\title{Positive characteristic Manin-Mumford theorem}
\author{Thomas Scanlon}
\thanks{Partially supported by NSF Grant DMS-0007190 and an
Alfred P. Sloan Fellowship}
\date{25 March 2003}
\email{scanlon@math.berkeley.edu}
\address{University of California, Berkeley \\
Department of Mathematics \\
Evans Hall \\
Berkeley, CA 94720-3480 \\
USA}
\maketitle

\begin{abstract}
We present the details of a model theoretic proof of 
an analogue of the Manin-Mumford conjecture for semiabelian 
varieties in positive characteristic. 
\end{abstract}

\section{Introduction}
The Manin-Mumford conjecture in its original form
(whose proof is originally due to Raynaud~\cite{Ray}) asserts 
that if $A$ is an abelian variety over a number field $k$
and $X \subseteq A$ is an irreducible subvariety of $A$, then 
$X(k^\alg)$ meets the torsion subgroup of $A(k^\alg)$ in a 
finite union of cosets of sugroups of the torsion group.  
If one replaces $k$ with a field of positive characteristic, then 
there are obvious counterexamples to the direct translation 
of this conjecture.  However, by isolating groups defined over 
finite fields appropriately, one can state and prove a positive
characteristic version of this conjecture. 

We should say a word or two about attributions for this theorem.  
The current author sketched the proof presented here in~\cite{ScBSL}.  
As the reader will see, given the dichotomy theorem for existentially 
closed difference fields~\cite{CHP} this proof follows very easily 
from Hrushovski's proof of the number field Manin-Mumford 
conjecture~\cite{HrMM}.  Pink and
Roessler gave an algebraic proof this theorem in~\cite{PR}.  While their
proof avoids appeals to the model theory of difference fields, it too uses 
some sophisticated arguements (involving, for instance, formal group 
arguments).  Pillay presented a very elementary proof of
the function field Mordell-Lang conjecture using an analysis of 
algebraic $D$-groups~\cite{PinML} and then  transposed this argument to the
context of algebraic $\sigma$-groups to reprove the Manin-Mumford conjecture
over number fields~\cite{PiMM}. 
The student working group supervised by Pillay and Scanlon
at the 2003 Arizona Winter School completed a very elementary proof of the 
main theorem of this note along the lines of Pillay's characteristic zero 
proof.  Some details of this argument are available in a streaming video on the
Southwestern Center's webpage.

\section{Statement of the main theorem}

In this section we state the main theorem of this note.  Before doing so
we recall a definition from~\cite{HrML}.

\begin{Def}
Let $K$ be an algebraically closed field of characteristic $p$.  
Let $G$ be a commutative algebraic group over $K$ and $X \subseteq G$
an irreducible subvariety.  We say that $X$ is \emph{special} if there 
are 
\begin{itemize}
\item $H \leq G$ an algebraic subgroup,
\item $H_0$ an algebraic group defined over $\FF_p^\alg$,
\item a point $a \in G(K)$, 
\item a subvariety $X_0 \subseteq G_0$ defined over $\FF_p^\alg$, and 
\item a morphism of algebraic groups $h:G \to G_0$
\end{itemize}

such that $X = a + h^{-1} X_0$.
\end{Def}

With the definition of \emph{special} in place we can state the 
positive characteristic version of the Manin-Mumford conjecture.

\begin{thm}
\label{MMp}
Let $K$ be an algebraically closed field of characteristic $p$.   
Let $G$ be a semiabelian variety over $K$ and $X \subseteq G$ a
closed subvariety.  Then the Zariski closure of $X(K) \cap G(K)_\tor$
is a finite union of special subvarieties.
\end{thm}

\section{The proof}

In this section we prove Theorem~\ref{MMp}.

\subsection{Difference equations for the torsion}
\begin{lem}
\label{prank}
There is a discrete valuation ring $R \subseteq K$
with a finite residue field $\FF_q$ and a semiabelian 
model $\fG$ of $G$ over $R$ for which the $p$-rank of the
special fibre of $\fG$ is equal to the $p$-rank of the 
generic fibre, $G$.
\end{lem}
\begin{proof}
Choose any finitely generated subring $S$ over which 
we have a semiabelian model $\fG$ of $G$.  Let $S' := S(G[p](K^\alg))$.  
Let $r$ be the $p$-rank of $G$ ($= \dim_{\FF_p} G[p](K^\alg)$). 
Let $\gamma_1, \ldots, \gamma_r \in G[p](K^\alg) = \fG[p](S')$ be a
basis for the (physical) $p$-torsion on $G$.   The set $U$ of 
primes $\fp \in {\rm Spec}(S')$ such that the image of 
$\gamma_1, \ldots, \gamma_r$ remain linearly independent in 
$(\fG \otimes S'/\fp)[p](S'/\fp)$ is
open in the Zariski topology. Take $\fm \in U$
any smooth closed point.  Take $R$ to be $S'_{\fm}$.
\end{proof}

\begin{lem}
\label{inttor}
Let $R$ be a discrete valuation ring of characteristic $p$
with residue field $\FF_q$ and field of quotients $K$.  
Let $S$ be the maximal unramified algebraic extension of $R$ and
let $S' := S^{p^{-\infty}} := \{ y \in K^\alg : (\exists n \in \NN) y^{p^n} \in S \}$
be the perfection of $S$.  
Then for any semiabelian scheme $G$ over $R$, the natural map 
$G(S')_{\rm tor} \to G(K^\alg)_\tor$ is an isomorphism.  
\end{lem}
\begin{proof}
Of course, this map is an injection.  So, we must show that it is a surjection.
For any finite \'{e}tale group scheme $F$ over $R$, Hensel's lemma shows that
$F(S) \hookrightarrow F(K^\alg)$ is an isomorphism.  For each $n \in \ZZ_+$,
Consider the connected-\'{e}tale sequence over $S'$: 

$$\begin{CD} 0 @>>> G[n]^0 @>>> G[n] @>>> G[n]_{\et} @>>> 0 \end{CD}$$

Over a perfect ring, this sequence splits and the group of rational points 
in a connected finite flat group scheme over a domain is trivial.  
Thus, $G[n](S') \cong G[n]_{\et}(S') \cong G[n]_{\et}(K^\alg) \cong G[n](K^\alg)$.
 
For each $n \in \NN$ the group scheme $G[n]$ is over $S$, thus if $\cO$ is the
integral closure of $R$ in $K^\alg$, we have $G[n](\cO) \cong G[n](K^\alg)$.  
As the torsion group is the direct limit of the $n$-torsion groups, we conclude
$G(S')_\tor \cong G(K^\alg)_\tor$.
\end{proof}

\begin{lem}
\label{toreqloc}
Let $R$ be a discrete valuation ring of characteristic $p$
with residue field $\FF_q$ and field of quotients $K$.
Let $G$ be a semiabelian scheme over $R$ for which the $p$-rank
of the generic fibre is equal to the $p$-rank of the 
special fibre. 
There is a polynomial $P(X) \in \ZZ[X]$ and an automorphism
$\sigma$ of $K^\alg$ fixing $K$ such that $P(\sigma)$ 
vanishes on $G(K^\alg)_\tor$ and no root of $P$ in $\CC$ 
is a root of unity. 
\end{lem}
\begin{proof}
On the special fibre $\overline{G}$ of $G$ the $q$-power Frobenius
induces an endomorphism $F:\overline{G} \to \overline{G}$.  
As such, the subring of ${\rm End}(\overline{G})$ generated by 
$F$ is a finite integral extension of $\ZZ$.  Let $P(X) \in \ZZ[X]$ 
be the minimal monic polynomial of $F$ over $\ZZ$. By the 
Weil conjectures for $\overline{G}$, no complex root of $P$
is a root of unity.  

The completion of $R$ is isomorphic to $\FF_q[[\epsilon]]$.  
Let $\rho:\FF_q^\alg[[\epsilon]] \to \FF_q^\alg[[\epsilon]]$ 
be defined by $$\sum_{i \geq 0} x_i \epsilon^i \to \sum_{i \geq 0} x_i^q \epsilon^i$$
Extend $\rho$ to
$\widetilde{\rho}:\FF_q^\alg((\epsilon))^\alg \to \FF_q^\alg((\epsilon))^\alg$ 
and let $\sigma := \widetilde{\rho} \upharpoonright_{K^\alg}$ be the restriction
of $\widetilde{\rho}$ to $K^\alg$.  

This choice of $P$ and $\sigma$ works.  By Lemma~\ref{inttor} every 
torsion point in $G(K^\alg)$ is integral over 
$S'$, the perfection of the maximal algebraic unramified extesion of $R$.
Our hypothesis on the $p$-rank 
implies that for each $n \in \ZZ_+$ the reduction map induces an isomorphism 
$G[n](S') \cong \overline{G}(\FF_q^\alg)$.  Moreover, as we have chosen $\sigma$
to lift $F$, if we regard $G(S')$ as a $\widehat{\ZZ}$-module with the generator
acting as $\sigma$ and $\overline{G}(\FF_q^\alg)$ as a $\widehat{\ZZ}$-module 
with the generator acting as $F$, then isomorphism 
$G(S')_\tor \to \overline{G}(\FF_q^\alg)$ is an isomorphism of $\widehat{\ZZ}$-modules.
$P$ is defined so that $P(F) \equiv 0$ on $\overline{G}(\FF_q^\alg)$.  Thus, 
$P(\sigma)$ vanishes on $G(S')_\tor = G(K^\alg)_\tor$.
\end{proof}

\subsection{Finite rank $\sigma$-algebraic groups}

\begin{lem}
\label{toreqglo}
Let $K = K^\alg$ be an algebraically closed field of characteristic $p > 0$ 
and $G$ a semiabelian variety over $K$.  There is a polynomial $P(X) \in \ZZ[X]$
having no roots of unity amongst its complex roots and an automorphism 
$\sigma:K \to K$ such that $G$ is defined over the fixed field of $\sigma$
and $P(\sigma)$ vanishes on $G(K)_\tor$.
\end{lem}

\begin{proof}
Giving $G$ a specific quasi-projective presentation, we find a finitely 
generated ring $R \subseteq K$ over which $G$ is semiabelian.  By Lemma~\ref{prank}
we may find a smooth closed point $\fm \subseteq R$ for which $G$ and $G_\fm$ 
have the same $p$-rank.  Applying Lemma~\ref{toreqloc} we obtain the requisite 
polynomial $P$ and automorphism $\sigma$.
\end{proof}

\begin{Def}
Let $K = K^\alg$ be an algebraically closed field and $k \leq K$ the 
algebraic closure of the prime field in $K$.  We say that the semiabelian
variety $G$ defined over $K$ is \emph{isotrivial} if there is a semiabelian
variety $G_0$ defined over $k$ and
a purely inseparable isogeny $\psi:G \to G_0$ defined over $K$.  (Equivalently, 
there is a purely inseparable isogeny $\vartheta:G_0 \to G$ defined
over $K$.)
\end{Def}

\begin{lem}
\label{isodes}
Let $M$ be an algebraically closed field and $K_1, K_2 \leq M$ algebraically
closed subfields which are algebraically independent over their intersection 
$K_1 \cap K_2$.  Suppose that $A$ is an algebraic group defined over 
$K_1$ and that there are an algebraic group $B$ defined over $K_2$ and 
a surjective map of algebraic groups $g:A \to B$ with $(\ker g)_{\rm red}$ defined 
over $K_1$ then there is an algebraic
group $B_0$ defined over $K_1 \cap K_2$ and a surjective morphism 
$h:A \to B_0$ defined over $K_1$ with $(\ker g)_{\rm red} = (\ker h)_{\rm red}$. 
\end{lem}
\begin{proof}
Choosing a presentation of $A$ over $K_1$, $B$ over $K_2$, and $g$ over $M$, 
we may express the assertion ``$g$ is a surjective map of algebraic groups 
from $A$ to $B$'' as a sentence in the language of fields with parameters from
$M$.  As the theory of algebraically closed fields is model complete, 
we may assume that all of the necessary parameters come from the algebraic
closure of the compositum of $K_1$ and $K_2$.  Separating the parameters 
and using quantifiers to speak about algebraic extensions, we may write
this sentence as $\varphi(a;b)$ where $a$ is a tuple from $K_1$, $b$ is a 
tuple from $K_2$, and $\varphi(x;y)$ is a formula of the language of rings 
having no extra parameters.  The formula $\varphi(x;y)$  asserts 
``$A_x$ is an algebraic group, $B_y$ is an algebraic group, there are parameters
$z$ satisfying a particular algebraic relation over $x$ and $y$ such that 
$g_{z}:A_x \to B_y$ is surjective and (on points) $\ker g_z = K_x$.''  

The formula $\varphi(x;y)$ is represented in ${\rm tp}(a/K_2)$, but $K_1$ and 
$K_2$ are free over $K_1 \cap K_2$.  Thus, $\varphi(x;y)$ is represented in 
${\rm tp}(a/K_1 \cap K_2)$.  That is, we can find a tuple $c$ from $K_1 \cap K_2$
for which $\varphi(a;c)$ holds.  This gives the result.
\end{proof}

\begin{lem}
\label{descent}
Let $(\U,+,\times,0,1,\sigma) \models \ACFA_p$ be an existentially 
closed difference field of characteristic $p > 0$.  Let $A$ be a semiabelian
variety defined over ${\rm Fix}(\sigma)$.  Suppose that there there are nonzero 
integers $m$ and $n$ such that $A$ is isogenous to a semiabelian variety 
defined over $\Fix (\sigma^n \tau^m)$, then $A$ is isogenous to a semiabelian variety 
defined over a finite field.
\end{lem}
\begin{proof}
The fields $\Fix(\sigma)$ and $\Fix(\sigma^n \tau^m)$ are orthogonal, and 
in particular, algebraically independent.  We have $\Fix(\sigma)^\alg =  \bigcup_{N \geq 0} 
\Fix(\sigma^N)$ and $\Fix(\sigma^n \tau^m)^\alg = \bigcup_{M \geq 0} \Fix(\sigma^nM \tau mM)$
so that $\Fix(\sigma)^\alg \cap \Fix(\sigma^n \tau^m)^\alg = \FF_p^\alg$.
Thus, using the fact that every algebraic subgroup of $A$ is defined over
${\rm Fix}(\sigma)^\alg$ by Lemma~\ref{isodes} we see that $A$ is isogenous to a
semiabelian variety defined over $\FF_p^\alg$.
\end{proof}

We recall also the definition of (quantifier-free) modularity in 
existentially closed difference fields.

\begin{Def}
Let $(\U,+,\times,0,1,\sigma) \models \ACFA$ be an existentially closed
difference field and $G$ an algebraic group over $\U$.  A subgroup 
$\Gamma \leq G(\U)$ is said to be \emph{modular} if for every natural 
number $n$ and every $\sigma$-algebraic subvariety $X \subseteq G^n$
the set $X(\U) \cap \Gamma^n$ is a finite union of cosets of subgroups of 
$\Gamma^n$.
\end{Def}

While it is not immediately clear from the definition, if $\Gamma$ and $\Xi \leq G(\U)$ 
are modular \emph{definable} groups, then so is the group generated by $\Gamma$
and $\Xi$. 

\begin{Def}
Let $(\U, +, \times, 0, 1, \sigma) \models \ACFA_p$ be an existentially 
closed difference field of characteristic $p > 0$. Denote the $p$-power
Frobenius map by $\tau:\U \to \U$.  Let $G$ be a commutative algebraic group
defined over $\U$ and $\Gamma \leq G(K)$ a definable subgroup.
We say that $\Gamma$ is \emph{essentially algebraic} if there is 
group $H$ of the form $\sum_{i=1}^m \psi_i(H_i(\Fix(\sigma^{n_i} \tau^{m_i}))$
where $n_i > 0$, $m_i \in \ZZ$, $H_i$ is an algebraic group over 
$\Fix(\sigma^{n_i}\tau^{m_i})$, and $\psi_i:H_i \to G$ is a map of 
algebraic groups with finite kernel such that $\Gamma/(H \cap \Gamma)$ is finite.

We say that $\Gamma$ is \emph{strongly essentially algebraic} if 
each of the $H_i$s may be taken to be defined over a finite field.
\end{Def}

\begin{lem}
\label{sea}
Let $(\U,+,\times,0,1,\sigma) \models \ACFA_p$ be an existentially closed
difference field of characteristic $p > 0$.  Let $G$ be a semisimple semiabelian 
variety over ${\rm Fix}(\sigma)$ and $Q(X) \in \ZZ[X]$ a polynomial with integral 
coefficients having no roots of unity amongst its complex roots.  
Let $\Gamma := \ker P(\sigma) \leq G(\U)$.  Then $\Gamma$ is an almost direct
sum of a modular group and a strongly essentially algebraic group.
\end{lem}

\begin{proof}
By hypothesis, we may 
write $G$ as an almost direct sum of simple semiabelian varieties $G = \sum S_i$
where each $S_i$ is either an abelian variety or an algebraic torus isomorphic
to $\Gm$.  As a general rule, the groups $S_i$ might not be 
invariant under $\sigma$, but they are always defined over a finite extension 
of the field of definition of $G$ and are therefore defined over ${\rm Fix}(\sigma^N)$
for some natural number $N$.  As there are only finitely many groups in this
decomposition, we may choose $N$ so that all of the $S_i$s an isomorphisms
$T_i \cong \Gm^{g_i}$ are defined over ${\rm Fix}(\sigma^N)$.

Factor $P(X)$ over $\CC$ as $P(X) = \prod_j (X - \alpha_j)$.  
Set $Q(X) := \prod_j (X - \alpha_j^N) \in \ZZ[X]$.  Then 
$\ker P(\sigma) \leq \ker Q(\sigma^N)$. If we manage to show that 
$\ker Q(\sigma^N) = M + E$ where $M$ is modular and $A$ is essentially 
algebraic, then $\ker P(\sigma) = (M \cap \Gamma) + (E \cap \Gamma)$
expresses $\Gamma$ as an almost direct sum of a modular group and
an essentially algebraic group.  Thus, replacing $\sigma$ with $\sigma^N$ and 
$P$ with $Q$ we may assume that each $S_i$ and each splitting $T_i \cong \Gm^{g_i}$
is already defined over $\Fix(\sigma)$.

Likewise, if we could show that each group $\ker P(\sigma) \upharpoonright S_i$
has the form $M_i + E_i$ where $M_i$ is modular and $E_i$ is essentially algebraic,
then we would conclude that $\ker P(\sigma) = (\sum M_i) + (\sum E_i)$ is the 
almost direct sum of the modular group $\sum M_i$ and the essentially algebraic
group $E_i$.  Thus, we may assume that $n = 1$ and that $G = S_1$.  We drop 
the subscript from now on.

One knows that every definable subgroup of
$\ker P(\sigma)$ in a simple semiabelian variety is commensurable with a group 
of the form $\ker R(\sigma)$ for some polynomial $R(X) \in {\rm End}(G)[X]$ over
the endomorphism ring of $G$ with $R$ dividing $P$~\cite{C}.  In particular, if one were to
factor $P(X) = \prod R_j(X)$ with each $R_j \in {\rm End}(G)[X]$ irreducible, then 
$\ker P(\sigma)$ is (up to finite index) the sum of the kernels $\ker R_j(\sigma)$.
Now, each of the groups $\ker R_j(\sigma)$ is minimal, so that the dichotomy theorem
implies that either $\ker R_j(\sigma)$ is modular or there is a minimal definable
field $k$, an algebraic group $H_0$ defined over $k$, and an isogeny $\psi:H_0 \to G$
such that $\psi(H(k)) \cap \ker R_j(\sigma)$ is infinite.    We consider 
two separate cases $k = \Fix(\sigma)$ or $k = \Fix(\sigma^n \tau^m)$ for $n > 0$ and 
$m \neq 0$.  

In the first case, we observe that $\psi$ is defined over $k^\alg$ and therefore
over $\Fix(\sigma^N)$ for some $N$.  The map $\psi$ takes $H(\Fix(\sigma^N))$
to a subgroup of finite index in $G(\Fix(\sigma^N))$.  
Factoring $R_j(X) = \prod (X - \beta_\ell)
\in \CC[X]$ and setting $Q_j = \prod (X - \beta_\ell^N) \in {\rm End}(G)[X]$,
we see that $\ker Q_j(\sigma^N)$ is commensurable with 
$\sum_{t = 0}^{N-1} \sigma^t \ker R_j(\sigma)$.   We conclude, that 
$\ker Q_j(\sigma^N)$ is commensurable with $G(\Fix(\sigma^N))$, but then 
some root of $Q_j$ (and therefore $R_j$ and also $P$) is a root
of unity contrary to our hypothesis.

In the second case, by Lemma~\ref{descent} (possibly at the cost of replacing 
$\sigma$ with a power) we see that may assume that $H$ 
is defined over a finite field.  Then by definition, $\ker R_j$ is strongly
essentially algebraic.
\end{proof}

We are in a position now to prove Theorem~\ref{MMp}.

\begin{proof}
Working by noetherian induction on $X$, we may assume that $X$ is irreducible and
that $X(K) \cap G(K)_\tor$ is Zariski dense in $X$.  Passing to the quotient 
by the stabilizer of $X$, we may assume that $X$ has a trivial stabilizer.  

Let $(\U,+,\times,0,1,\sigma) \models \ACFA$ be an existentially closed difference
field with $K \leq \U$, $\sigma(G) = G$, and $P(X) \in \ZZ[X]$ a polynomial over
the integers with no roots of unity amongst its complex roots and $P(\sigma)$
vanishing on $G(K)_\tor$.  Let $E \leq \ker P(\sigma) =: \Gamma$ be an essentially algebraic
subgroup of maximal dimension.  By the socle theorem, $X(K) \cap \Gamma$ is
contained in finitely many cosets of $E$.  Translating, we may assume that $X$ is a 
subvariety of the Zariski closure of $E$.  By Lemma~\ref{sea} $E$ is strongly
essentially algebraic so that there is some algbraic group $H_0$ over a finite
field and an isogeny $\psi:E \to H_0$ which we may assume to be purely 
inseparable. Let $X_0$ be the Zariski closure of $\psi(X)$.  As $\psi$ takes
torsion points to torsion points and $H_0(\U)_\tor = H_0(\FF_p^\alg)$, 
we see that $X_0$ is defined over $\FF_p^\alg$.  Thus, $X$ is special.
\end{proof}

\end{document}